\documentclass{article} % For LaTeX2e
\usepackage{authblk}
\usepackage{hyperref}
\usepackage{url}

\usepackage{amsmath,amsfonts,amssymb}
\usepackage{mathrsfs}

\newcommand{\set}[1]{\left\{ #1 \right\}}
\newcommand{\norm}[1]{\left\Vert #1 \right\Vert}

\newcommand{\tr}[1]{\mathrm{Tr} \left( #1 \right)}
\newcommand{\abs}[1]{\left\vert #1 \right\vert}

\newcommand{\ttrue}{\theta^\natural}
\newcommand{\that}{\hat{\theta}}
\newcommand{\zhat}{\hat{z}}
\newcommand{\ztrue}{z^\natural}

\newcommand{\fset}{\mathcal{F}_g ( \ttrue )}
\newcommand{\fcone}{\overline{\fset}}

\usepackage{amsthm}
\theoremstyle{plain}% default
\newtheorem{thm}{Theorem}[section]
\newtheorem{lem}[thm]{Lemma}
\newtheorem{prop}[thm]{Proposition}
\newtheorem{cor}[thm]{Corollary}

\theoremstyle{definition}
\newtheorem{defn}{Definition}[section]

\theoremstyle{remark}
\newtheorem*{rem}{Remark}

\usepackage{xcolor}

\title{A Geometric View on Constrained $M$-Estimators}

\author{Yen-Huan~Li, Ya-Ping~Hsieh, Nissim Zerbib and Volkan~Cevher}

\affil{Laboratory for Information and Inference Systems\\
\'{E}cole Polytechnique F\'{e}d\'{e}rale de Lausanne}

\date{\empty}

\begin{document}

\maketitle

\begin{abstract} 
We study the estimation error of constrained $M$-estimators, and derive explicit upper bounds on the expected estimation error determined by the Gaussian width of the constraint set. Both of the cases where the true parameter is on the boundary of the constraint set (matched constraint), and where the true parameter is strictly in the constraint set (mismatched constraint) are considered. For both cases, we derive novel universal estimation error bounds for regression in a generalized linear model with the canonical link function. Our error bound for the mismatched constraint case is minimax optimal in terms of its dependence on the sample size, for Gaussian linear regression by the Lasso.
\end{abstract} 

\section{Introduction} \label{sec_formulation}

Consider a general statistical estimation problem. Let $( y_1, \ldots, y_n )$ be a sample following a probability distribution $\mathbb{P}_{\ttrue}$ in a given class $\mathcal{P} := \set{ \mathbb{P}_{\theta} : \theta \in \mathbb{R}^p }$. We are interested in estimating the parameter $\ttrue$, given $( y_1, \ldots, y_n )$ and $\mathcal{P}$, under the high-dimensional setting where $n < p$.

If $\ttrue$ is known to satisfy $g ( \ttrue ) \leq c$ for some continuous convex function $g$ and positive constant $c$, we can consider a constrained $M$-estimator of the form
\begin{equation}
\that \in \arg \min_{ \theta } \set{ f_n ( \theta ) : \theta \in \mathcal{G} }, \quad \mathcal{G} := \set{ \theta \in \mathbb{R}^p : g ( \theta ) \leq c }. \label{eq_that}
\end{equation}
We assume that $f_n$ is a continuously differentiable convex function, and the constraint set $\mathcal{G}$ is non-empty. For example, the Lasso \cite{Tibshirani1996} corresponds to
\begin{equation}
f_n ( \theta ) := \frac{1}{2 n} \sum_{i = 1}^n \left( y_i - \left\langle a_i, \theta \right\rangle \right)^2, \quad \mathcal{G} := \set{ \norm{ \theta }_1 \leq c }, \label{eq_LS} 
\end{equation}
for some $a_1, \ldots, a_n \in \mathbb{R}^p$ and positive constant $c$
%\footnote{The $\ell_1$-regularized least squares estimator is sometimes also called the Lasso, while the original Lasso defined in \cite{Tibshirani1996} is of the constrained form.}
. A matrix $\Theta \in \mathbb{R}^{d \times d}$ can be vectorized as a corresponding vector $\theta \in \mathbb{R}^p$, $d^2 = p$. In the low-rank matrix recovery problem \cite{Candes2011b,Gunasekar2014}, a popular estimator corresponds to
\begin{equation}
f_n ( \Theta ) := \frac{1}{2 n} \sum_{i = 1}^n \left( y_i - \tr{ A_i^T \Theta } \right)^2, \quad \mathcal{G} := \set{ \norm{ \Theta }_* \leq c }, \label{eq_matrix_lasso}
\end{equation}
for some $A_1, \ldots, A_n \in \mathbb{R}^{d \times d}$ and positive constant $c$, where $\norm{ \cdot }_*$ denotes the nuclear norm. In general, $f_n$ can be the normalized negative log-likelihood function, or any properly defined function, and $g$ depends on the \textit{a priori} information on the structure of the parameter $\ttrue$ \cite{Bach2013a,Chandrasekaran2012,ElHalabi2015}.

One can also consider a penalized $M$-estimator, given by
\begin{equation}
\that_{\text{penalized}} \in \arg \min_{ \theta \in \mathbb{R}^p } \set{ f_n ( \theta ) + \rho_n g ( \theta ) }, \label{eq_penalized_est}
\end{equation}
for some positive constant $\rho_n$. The penalized $M$-estimator can be computed by fast proximal methods, provided that the proximal mapping of $g$ is easy to compute \cite{Beck2009,Nesterov2013}. This condition, however, is not always satisfied. For example, if $g$ is the nuclear norm, computing the corresponding proximal mapping requires a full singular value decomposition (SVD) in the first few iterations, and hence is not scalable with the parameter dimension. In contrast, if we consider a constrained $M$-estimator and compute it by the Frank-Wolfe algorithm, each iteration of the algorithm requires a linear minimization oracle (LMO), which %only needs the top eigenvalue of a matrix, and 
can be approximated efficiently by Lanczos' algorithm \cite{Jaggi2013}. The paper \cite{Zhang2013} also shows that when $g$ is a structured sparsity regularizer, the LMO can be much easier to compute than the proximal mapping. 
%Despite the computational superiority of constrained $M$-estimators in some applications, current results about the performance of a constrained $M$-estimator is quite limited. We shall address this in Section \ref{sec_related_work}.

If we consider a constrained $M$-estimator, setting the value of the constant $c$ in (\ref{eq_that}) becomes a practical issue. 
For the case $c < g ( \ttrue )$, the estimation error is obviously bounded below by the distance between $\ttrue$ and the constraint set $\mathcal{G}$, and hence estimation consistency is impossible.
Ideally we would like to set $c = g ( \ttrue )$, while in practice $g ( \ttrue )$ is seldom known. The last case is when we have some estimate on $g( \ttrue )$, and choose $c$ such that $c > g ( \ttrue )$. Some natural questions arise: Is estimation consistency possible? How fast will the estimation error decay with the sample size $n$? Does setting $c > g ( \ttrue )$ result in larger estimation error than setting $c = g ( \ttrue )$? We review related works in Section \ref{sec_related_work}, which shows that answers existed only for specific cases even when $c = g ( \ttrue )$.

In this paper, we provide a unified analysis for constrained $M$-estimators.
%, in a manner \cite{Negahban2012} does for penalized estimators
Specifically, 

\begin{itemize}
\item We propose an elementary framework for analyzing any $M$-estimator applied to any statistical model in Section \ref{sec_framework}.
\item We obtain universal error bounds in terms of the Gaussian width, valid \emph{for all} canonical GLMs. We consider the matched constraint case ($c = g ( \ttrue )$) in Section \ref{sec_matched}, and the mismatched constraint case ($c > g ( \ttrue )$) in Section \ref{sec_mismatch}.
\item To illustrate the universal error bounds, we specialize the universal error bound to Gaussian linear regression with arbitrary convex constraint, and regression in canonical GLMs with the $\ell_1$-constraint in Section \ref{sec_appl}, and obtain explicit results. 
\item Our error bound for the Lasso applied to the Gaussian linear model is optimal in the minimax sense (cf. Section \ref{sec_mismatch_further}).
% \item Our results are expressed in fundamental geometric quantities of the constraint set, such as the Gaussian squared-complexity and the Gaussian width.
\end{itemize}

Existing results for penalized $M$-estimators \cite{Banerjee2015,Bickel2009,Buhlmann2011,Honorio2014,Kakade2010,Negahban2012,Geer2013}, which are for deterministic $\rho_n$'s, cannot directly recover our results, and vice versa.
Indeed, by Lagrange duality, there exists some $\rho_n > 0$ such that the constrained $M$-estimator in (\ref{eq_that}) is equivalent to the penalized $M$-estimator in (\ref{eq_penalized_est}). 
This correspondence, however, holds \emph{only for given realization of the sample $( y_1, \ldots, y_n )$}, and hence $\rho_n$ is a random variable depending on the sample. 
Conversely, for any penalized $M$-estimator $\hat{\theta}_{\text{penalized}}$ for some $\rho_n > 0$, there exists a constant $c = g ( \hat{\theta}_{\text{penalized}} )$ such that the corresponding constrained $M$-estimator (\ref{eq_penalized_est}) is equivalent to $\hat{\theta}_{\text{penalized}}$. 
Note that $c = g ( \hat{\theta}_{\text{penalized}} )$ is again a random variable and dependent on the sample. 
We are not aware of any existing work on characterizing the correspondence between the two formulations.

\section{Related Works} \label{sec_related_work}

% The estimation errors of constrained $M$-estimators under the high-dimensional setting had been studied in various ways \cite{Koltchinskii2013a,Mendelson2014,Oymak2013a,Oymak2013,Plan2013a,Plan2015,Plan2014a,Vershynin2014}. We briefly address the limitations of the cited works below.

% Some of the existing works focus on the fundamental linear model.
% The work \cite{Chatterjee2014} considers the constrained least squares estimator for estimating the mean of a Gaussian process, and provides a concentration inequality for the estimation error. Note that this scenario corresponds to the case where $n = p$, and is not actually high-dimensional.
In \cite{Oymak2013a,Oymak2013}, the authors derived sharp estimation error bounds for regression in the linear model by constrained least squares (LS) estimators. 
The analysis in \cite{Vershynin2014} provides a minimax estimation error bound for the same setting
%\footnote{The paper \cite{Vershynin2014} also reviews the result in \cite{Plan2014a} for possibly non-linear statistical models, which we shall address in the next paragraph.}
.
There are some related works on learning a function in a function class \cite{Koltchinskii2013a,Mendelson2014}. When the function class is linearly parametrized by vectors in $\mathbb{R}^p$, and the function corresponding to $\theta^\natural$ is in the function class, the $L_2$-estimation error in the function class may be translated into the $\ell_2$-estimation error with respect to $\theta^\natural$.
A common limitation of \cite{Koltchinskii2013a,Mendelson2014,Oymak2013,Oymak2013a,Vershynin2014} is that the results are not extendable to general non-linear statistical models. 

Another research direction considers constrained estimation in possibly non-linear statistical models \cite{Plan2013a,Plan2015,Plan2014a}.
A constrained $M$-estimator for logistic regression was proposed and analyzed in \cite{Plan2013a}. 
In \cite{Plan2014a}, the authors proposed and analyzed a universal projection-based estimator for regression in generalized linear models (GLMs). 
In \cite{Plan2015}, the authors analyzed the performance of the constrained LS estimator in GLMs.
A common limitation of \cite{Plan2013a,Plan2015,Plan2014a} is that the results are valid only for the specific proposed estimators, and they do not even apply to the constrained maximum-likelihood (ML) estimator, which is the most popular approach in practice. Moreover, the proposed estimators in \cite{Plan2013a,Plan2015,Plan2014a} can only recover the true parameter up to a scale ambiguity.

% We would like to emphasize a technical issue in the cited works, regarding whether the constraint is matched or not. 
We say that the constraint is \emph{matched} if $\theta^\natural$ lies on the boundary of $\mathcal{G}$ in (\ref{eq_that}) (or $c = g ( \theta^\natural )$), and \emph{mismatched} if $\ttrue$ lies strictly in $\mathcal{G}$ (or $c < g ( \ttrue )$). The analyses in \cite{Oymak2013a,Oymak2013} require the constraint to be matched, while in practice the exact value of $g( \theta^\natural )$ is seldom known. The constraint in \cite{Koltchinskii2013a} is always matched due to the special structure of quantum density operators. The error bounds in \cite{Plan2013a,Vershynin2014} can be overly pessimistic, because they hold for all $\theta^\natural \in \mathcal{G}$. The results in \cite{Mendelson2014,Plan2015,Plan2014a} do not require a matched constraint and depend on $\theta^\natural$; our result is of this kind. Recall that, however, \cite{Mendelson2014} is limited to specific statistical models, and \cite{Plan2015,Plan2014a} are limited to specific $M$-estimators.

\section{A Geometric Framework} \label{sec_framework}

\subsection{Basic Idea} \label{sec_basic_idea}

To illustrate the basic idea of our framework, let us start with a simple setting, where $f_n$ is strongly convex with parameter $\mu > 0$, i.e., 
\begin{equation}
\left\langle \nabla f_n ( y ) - \nabla f_n ( x ), y - x \right\rangle \geq \mu \norm{ y - x }_2^2, \notag
\end{equation} 
for any $x, y \in \mathrm{dom}\, f$. Note that then $\that$ is uniquely defined.

Define $\iota_g: \mathbb{R}^p \to \mathbb{R} \cup \set{ + \infty }$ as the indicator function of the constraint set $\mathcal{G}$; that is, $\iota_{\mathcal{G}} ( \theta ) = 0$ if $\theta \in \mathcal{G}$, and $\iota_{ \mathcal{G} } ( \theta ) = + \infty$ otherwise.
By the strong convexity of $f_n$, we have
\begin{equation}
\left\langle \nabla f_n ( \that ) - \nabla f_n ( \ttrue ), \that - \ttrue \right\rangle \geq \mu \norm{ \that - \ttrue }_2^2. \label{eq_strong_convexity}
\end{equation}
By the convexity of $\iota_g$, or the monotonicity of the subdifferential mapping, we have
\begin{equation}
\left\langle \zhat - \ztrue, \that - \ttrue \right\rangle \geq 0, \label{eq_monotonicity}
\end{equation}
for any $\zhat \in \partial \iota_g ( \that )$, and any $\ztrue \in \partial \iota_g ( \ttrue )$. Summing up (\ref{eq_strong_convexity}) and (\ref{eq_monotonicity}), we obtain
\begin{equation}
\left\langle \nabla f_n ( \that ) + \zhat - \nabla f_n ( \ttrue ) - \ztrue, \that - \ttrue \right\rangle \geq \mu \norm{ \that - \ttrue }_2^2, \notag
\end{equation}
for any $\zhat \in \partial \iota_g ( \that )$. By the optimality condition of $\that$, there exists some $\hat{z} \in \partial \iota_{\mathcal{G}} ( \hat{\theta} )$ such that
\begin{equation}
0 = \nabla f_n ( \that ) + \hat{z}, \label{eq_optimality}
\end{equation}
and hence we have
\begin{equation}
\left\langle - \nabla f_n ( \ttrue ) - \ztrue, \that - \ttrue \right\rangle \geq \mu \norm{ \that - \ttrue }_2^2, \notag
\end{equation}
for any $\ztrue \in \partial \iota_g ( \ttrue )$. Since $\partial \iota_g ( \ttrue )$ is always a closed convex cone, we can choose $\ztrue = 0$ and obtain
\begin{equation}
\left\langle - \nabla f_n ( \ttrue )  , \that - \ttrue \right\rangle \geq \mu \norm{ \that - \ttrue }_2^2. \label{eq_bad_lower}
\end{equation}
Applying the Cauchy-Schwarz inequality to the left-hand side, we obtain
\begin{align}
\norm{ \nabla f_n ( \ttrue ) }_2 \norm{ \that - \ttrue }_2 \geq  \mu \norm{ \that - \ttrue }_2^2, \notag
\end{align}
or
\begin{equation}
\norm{ \that - \ttrue }_2 \leq \frac{1}{\mu} \norm{ \nabla f_n ( \ttrue ) }_2. 
\end{equation}
Taking expectations on both sides, we immediately obtain the following estimation error bound: 
\begin{equation}
\mathbb{E}\, \norm{ \that - \ttrue }_2 \leq \frac{1}{\mu} \mathbb{E}\, \norm{ \nabla f_n ( \ttrue ) }_2. \label{eq_bad_bound}
\end{equation}
The gradient at the true parameter, $\nabla f_n ( \ttrue )$, usually concentrates around $0$ with high probability. 

The simple error bound (\ref{eq_bad_bound}) is not desirable for two reasons: 
\begin{enumerate}
\item In the high-dimensional setting where $n < p$, $f_n$ cannot be strongly convex even for the basic LS estimator.
\item It does not depend on the choice of $g$.
\end{enumerate}

We address the first issue in Section \ref{sec_RSC}, and the second issue in Section \ref{sec_refined_bound}. 

\subsection{Restricted Strong Convexity} \label{sec_RSC}

Note that in order to facilitate the arguments in the previous sub-section, we only require (\ref{eq_strong_convexity}) to hold for $\that$ and $\ttrue$, instead of any two vectors in $\mathbb{R}^p$. Therefore, we only need $f_n$ to satisfy some \emph{restricted} notion of strong convexity. Similar (but not exactly the same) ideas had appeared in \cite{Chandrasekaran2012,Negahban2012}, and can be traced back to \cite{Bickel2009,Geer2007}.

\begin{defn}[Feasible Set and Feasible Cone]
The \emph{feasible set} of $g$ at $\ttrue$, denoted by $\mathcal{F}_g ( \ttrue )$, is given by
\begin{equation}
\mathcal{F}_g ( \ttrue ) := \mathcal{G} - \ttrue = \set{ \theta - \ttrue: \theta \in \mathcal{G} }. \notag
\end{equation}
The \emph{feasible cone} of $g$ at $\ttrue$, denoted by $\overline{\mathcal{F}_g ( \ttrue )}$, is defined as the conic hull of $\mathcal{F}_g ( \ttrue )$.
\end{defn}

By the definition of $\that$, the estimation error must satisfy $\that - \ttrue \in \mathcal{F}_g ( \ttrue )$.

\begin{defn}[Restricted Strong Convexity] \label{def_RSC}
The function $f_n$ satisfies the restricted strong convexity (RSC) condition with parameter $\mu > 0$ if
\begin{equation}
\left\langle \nabla f_n ( \ttrue + e ) - \nabla f_n ( \ttrue ), e \right\rangle \geq \mu \norm{ e }_2^2, \label{eq_RSC}
\end{equation}
for any $e \in \mathcal{F}_g ( \ttrue )$.
\end{defn}

If $f_n$ is twice continuously differentiable, we have a sufficient condition.

\begin{prop} \label{prop_RSC_hessian}
The function $f_n$ satisfies the RSC condition with parameter $\mu > 0$ if
\begin{equation}
\left\langle e, \nabla^2 f_n ( \ttrue + \lambda e ) e \right\rangle \geq \mu \norm{e}_2^2, \notag
\end{equation}
for all $\lambda \in [0,1]$ and all $e \in \fset$.
\end{prop}

The uniqueness of $\that$ and the derivation of the error bound in Section \ref{sec_basic_idea} are still valid even when $n < p$, as long as $f_n$ satisfies the RSC condition with some parameter $\mu > 0$.

\subsection{Refined Error Bound} \label{sec_refined_bound}

We address the dependence of the estimation error on the choice of $g$, and derive a refined error bound in this sub-section.

We note that 
\begin{equation}
\left\langle - \nabla f_n ( \ttrue ), \that - \ttrue \right\rangle = \left\Vert \Pi_{ \overline{ \that - \ttrue } } \left( - \nabla f_n ( \ttrue ) \right) \right\Vert_2 \norm{ \that - \ttrue }_2, \notag
\end{equation}
where $\Pi_{ \overline{ \that - \ttrue } } ( \cdot )$ denotes the projection onto the conic hull of $\set{\that - \ttrue}$ (which is a half-line or $\{ 0 \}$). This implies, by (\ref{eq_bad_lower}), 
\begin{equation}
\norm{ \Pi_{ \overline{ \that - \ttrue } } \left( - \nabla f_n ( \ttrue ) \right) }_2 \geq \mu \norm{ \that - \ttrue }_2. \notag
\end{equation}
The left-hand side, however, is not tractable due to its dependence on $\that$. As $\that - \ttrue \in \fcone$ by definition, we consider a looser bound: 
\begin{equation}
\norm{ \Pi_{\overline{\mathcal{F}_g ( \ttrue )}} ( - \nabla f_n ( \ttrue ) ) }_2 \geq \mu \norm{ \that - \ttrue }_2, \label{eq_concentration}
\end{equation}
where $\Pi_{\fcone} ( \cdot )$ denotes projection onto the feasible cone $\fcone$.

%The left-hand side, however, is not tractable due to its dependence on $\that$, so we consider a looser bound:
%\begin{equation}
%\sup_{ e \in \fset } \set{ \norm{ \Pi_{\overline{e}} \left( - \nabla f_n ( \ttrue ) \right) }_2 } \geq \mu \norm{ \that - \ttrue }_2, \notag
%\end{equation}
%where $\Pi_{\overline{e}} ( \cdot )$ denotes the projection onto the conic hull of $\set{e}$. 
%The left-hand side is simply the projection of $- \nabla f_n ( \ttrue )$ onto the feasible cone $\fcone$. 

Taking expectations on both sides, we obtain the following lemma.

%Now we address the dependence of the estimation error on the choice of $g$, based on Moreau's decomposition.
%
%\begin{thm}[Moreau's Decomposition \cite{Moreau1965}] \label{thm_Moreau}
%Let $\mathcal{K} \subseteq \mathbb{R}^p$ be a closed convex cone. Any vector $v \in \mathbb{R}^p$ can be uniquely decomposed as
%\begin{equation}
%v = \Pi_{\mathcal{K}} ( v ) + \Pi_{\mathcal{K}^\circ} ( v ), \notag
%\end{equation}
%where $\Pi_{\mathcal{K}} ( \cdot )$ and $\Pi_{\mathcal{K}^\circ} ( \cdot )$ denote the projection onto $\mathcal{K}$ and its polar cone $\mathcal{K}^\circ$, respectively. Furthermore, we have
%\begin{equation}
%\left\langle \Pi_{\mathcal{K}} ( v ), \Pi_{\mathcal{K}^\circ} ( v ) \right\rangle = 0. \notag
%\end{equation}
%\end{thm}
%
%By Moreau's decomposition and the definition of the polar cone, we have
%\begin{equation}
%\left\langle \Pi_{\overline{\mathcal{F}_g ( \ttrue )}} ( - \nabla f_n ( \ttrue ) ), \that - \ttrue \right\rangle \geq \left\langle - \nabla f_n ( \ttrue )  , \that - \ttrue \right\rangle. \notag
%\end{equation}
%Similarly as in Section \ref{sec_basic_idea}, by applying the Cauchy-Schwarz inequality on the left-hand side, we obtain the following estimation error bound: 
%\begin{equation}
%\norm{ \that - \ttrue }_2 \leq \frac{1}{\mu} \norm{ \Pi_{\overline{\mathcal{F}_g ( \ttrue )}} ( - \nabla f_n ( \ttrue ) ) }_2, \notag
%\end{equation}
%which implies the following lemma.

\begin{lem} \label{lem_fundamental}
Assume that $f_n$ satisfies the RSC condition with parameter $\mu > 0$. Then $\that$ is uniquely defined, and satisfies
\begin{equation}
\mathbb{E}\, \norm{ \that - \ttrue }_2 \leq \frac{1}{\mu} \mathbb{E}\, \norm{ \Pi_{\overline{\mathcal{F}_g ( \ttrue )}} ( - \nabla f_n ( \ttrue ) ) }_2. \notag
\end{equation}
\end{lem}

Since $- \nabla f_n ( \ttrue )$ is a descent direction of $f_n$, if its direction is coherent with the feasible cone $\fcone$, we may find some point $\that'$ far away from $\ttrue$ in the feasible set $\fset $ such that $f_n ( \that' )$ is much smaller than $f_n ( \ttrue )$, and hence the estimation error can be large. This provides an intuitive interpretation of the lemma.

Since projection onto a closed convex set is a non-expansive mapping, we have
\begin{equation}
\norm{ \Pi_{\overline{\mathcal{F}_g ( \ttrue )}} ( - \nabla f_n ( \ttrue ) ) }_2 \leq \norm{ \nabla f_n ( \ttrue ) }_2, \notag
\end{equation}
so the error bound is always no larger than the one in Section \ref{sec_basic_idea}.

Lemma \ref{lem_fundamental} is the theoretical foundation of the rest of this paper. 

\section{Estimation Error Bound in Terms of the Gaussian Width} \label{sec_matched}

We apply Lemma \ref{lem_fundamental} to constrained ML estimators in a GLM with the canonical link function. Examples of a canonical GLM include the Gaussian linear, logistic, gamma, and Poisson regression models.

Let $\ttrue \in \mathbb{R}^p$ be the parameter to be estimated, or the unknown vector of regression coefficients. In a canonical GLM, the negative log-likelihood of a sample $y$, given $\ttrue$, is of the form (up to scaling and shifting by some constants)
\begin{equation}
L ( y; \ttrue ) = y \left\langle a_i, \ttrue \right\rangle - b ( \left\langle a_i, \ttrue \right\rangle ), \notag
\end{equation}
where $a_1, \ldots, a_n \in \mathbb{R}^p$ are given, and we assume that $b$ is some given concave function. Let $( y_1, \ldots, y_n ) \in \mathbb{R}^n$ be the sample. The constrained ML estimator is given by (\ref{eq_that}) with 
\begin{equation}
f_n ( \theta ) := \frac{1}{n} \sum_{i = 1}^n L ( y_i, \theta ),  \label{eq_fn_glm}
\end{equation}
and $g$ being some continuous convex function. For simplicity, we consider the case where $c = g ( \ttrue )$ in this section; we address the case where $c > g ( \ttrue )$ in Section \ref{sec_mismatch}.

We specialize Lemma \ref{lem_fundamental} to the canonical GLM and obtain the following theorem.

\begin{defn}[Gaussian width \cite{Chandrasekaran2012,Mendelson2007,Tropp2014}] \label{def_gwidth}
Let $\mathcal{C} \subseteq \mathbb{R}^p$. The \emph{Gaussian width} of $\mathcal{C}$ is given by
\begin{equation}
\omega_t ( \mathcal{C} ) := \mathbb{E}\, \sup_{ v \in \mathcal{C} \cap t \mathcal{S}^{p - 1} } \set{ \left\langle h, v \right\rangle }, \notag
\end{equation}
where $h := ( h_1, \ldots, h_p )$ is a vector of i.i.d. standard Gaussian random variables, and $\mathcal{S}^{p-1}$ denotes the unit $\ell_2$-sphere in $\mathbb{R}^p$.
\end{defn}

\begin{thm} \label{thm_no_mismatch}
Consider the canonical GLM and the corresponding ML estimator described above for $c = g ( \ttrue )$. 
Assume that the entries of $a_1, \ldots, a_n$ are either all i.i.d. standard Gaussian or all i.i.d. Rademacher random variables (random variables taking values in $\set{ +1, -1}$ with equal probability), and $f_n$ satisfies the RSC condition for $\mu > 0$ with probability at least $1/2$. 
Then
\begin{equation}
\mathbb{E}\, \norm{ \hat{\theta} - \ttrue }_2 \leq 2 \sqrt{ 2 \pi } \, \sigma_{\max} \frac{\omega_1 ( \fcone )}{ \mu \sqrt{n}}, \notag
\end{equation}
where $\sigma_{\max} := \max_i \sqrt{ \mathrm{var}\, y_i }$.
\end{thm}

\begin{rem}
Note that the expectation is with respect to $A$ and $\varepsilon$, conditioned on the event that the RSC condition holds.
\end{rem}

%\begin{rem}
%The constant $2$ in the error bound is arbitrary. A finer choice of the constant depends on the probability that the RSC holds (c.f. Section \secProofNoMismatch).
%\end{rem}

The feasible cone $\overline{ \mathcal{F}_g ( \ttrue ) }$ coincides with the tangent cone of $g$ at $\ttrue$ defined in \cite{Chandrasekaran2012}. Therefore, to evaluate the estimation error bound, we only need to evaluate the Gaussian width of the corresponding tangent cone. We note that there are already many results for a variety of commonly used regularization functions, such as the $\ell_1$-norm, nuclear norm, total variation semi-norm, and general atomic norms \cite{Cai2013,Chandrasekaran2012,Foygel2014,Plan2013a,Rao2012,Vershynin2014}. 
Therefore, for most of the applications, we only need to \emph{plug in} an existing bound on the Gaussian width.

Finally, we would like to emphasize that the Gaussian width in Theorem \ref{thm_no_mismatch} comes from bounding the random process induced by the random gradient $\nabla f_n ( \theta^\natural )$ (cf. the proof of Theorem \ref{thm_no_mismatch}), instead of being a consequence of applying Gordon's Lemma. That is, our result is essentially different from those in \cite{Chandrasekaran2012,Oymak2013a,Oymak2013}.

\section{Effect of a Mismatched Constraint} \label{sec_mismatch}

In this section, we discuss the effect of a mismatched constraint for ML regression in a canonical GLM. Recall that the constraint set $\mathcal{G}$ is called \emph{mismatched} if $c > g ( \ttrue )$ in (\ref{eq_that}).

The notion of the RSC in Definition \ref{def_RSC} is no longer meaningful when the constraint set is mismatched. Take ML regression in the Gaussian linear model for example, for which the corresponding $f_n$ is given in (\ref{eq_LS}). Let $A \in \mathbb{R}^{n \times p}$ be defined as in Theorem \ref{thm_no_mismatch}. 
%By direct calculations, we get
%\begin{equation}
%\nabla f_n ( \theta ) = - \frac{1}{n} A^T ( y - A \theta ). \notag
%\end{equation}
The RSC condition requires
\begin{equation}
\left\langle \nabla f_n ( \ttrue + e ) - \nabla f_n ( \ttrue ), e \right\rangle = \frac{1}{n} \norm{A e}_2^2 \geq \mu \norm{ e }_2^2, \notag
\end{equation}
for some $\mu > 0$ and all $e \in \fcone$, where we say $e \in \fcone$ instead of $e \in \fset$ because $A$ is a linear operator. Since when the constraint is mismatched, $\fcone$ is the whole space $\mathbb{R}^p$, the RSC condition requires $A$ to be a non-singular matrix. This cannot be true in the high-dimensional setting, where $A \in \mathbb{R}^{n \times p}$ and $n < p$.

\textbf{Our Approach: } Let $t > 0$ and denote by $\mathcal{B}$ the unit $\ell_2$-ball in $\mathbb{R}^p$. We partition the feasible set $\fset$ as
\begin{equation}
\fset = ( \fset \cap t \mathcal{B} ) \cup ( \fset \setminus t \mathcal{B} ). \notag
\end{equation}
When $t$ is large enough, the conic hull of $( \fset \setminus t \mathcal{B} )$ will not be the whole space $\mathbb{R}^p$, so it is possible to have restricted strong convexity on $( \fset \setminus t \mathcal{B} )$ when $n < p$. If the error vector $\that - \ttrue$ lies in $( \fset \setminus t \mathcal{B} )$, we can obtain an error bound, say, $\tilde{t}$, as in Section \ref{sec_matched}; otherwise, if the error vector lies in $\mathcal{B}_t$, a na\"{\i}ve error bound is the radius of the ball, i.e., $t$. Finally, we can bound the estimation error from above by the maximum of $\tilde{t}$ and $t$. Note that $\tilde{t}$ is implicitly dependent on $t$.

The arguments in the previous paragraph can be made precise as in Lemma \ref{lem_mismatched}, which is an analogue of Lemma \ref{lem_fundamental} in the mismatched case. Lemma \ref{lem_mismatched} holds for arbitrary constrained $M$-estimators of the form (\ref{eq_that}) and statistical models.

\begin{lem} \label{lem_mismatched}
Suppose that for some $t > 0$, we have
\begin{gather}
\left\langle \nabla f_n ( \ttrue + e ) - \nabla f_n ( \ttrue ), e \right\rangle \geq \mu \norm{e}_2^2, \label{eq_RRSC}
\end{gather}
for some $\mu > 0$ and all $e \in \fset \setminus t \mathcal{B}$. Then 
\begin{equation}
\mathbb{E}\, \norm{ \that - \ttrue }_2 \leq t + \mathbb{E}\, \norm{ \Pi_{ \overline{ \fset \setminus t \mathcal{B} } } \left( - \nabla f_n ( \ttrue ) \right) }_2. \notag
\end{equation}
\end{lem}

We can also prove an analogue of Theorem \ref{thm_no_mismatch} for constrained ML regression in a canonical GLM.

\begin{cor} \label{cor_mismatched}
Consider the canonical GLM and the corresponding ML estimator described in Section \ref{sec_matched}, for $c > g ( \ttrue )$. Let $A$ be defined as in Theorem \ref{thm_no_mismatch} and let $t > 0$. Suppose that (\ref{eq_RRSC}) holds true with for some $\mu > 0$ with probability at least $1/2$.
Then we have
\begin{equation}
\mathbb{E}\, \norm{ \that - \ttrue }_2 \leq t +  2 \sqrt{ 2 \pi } \, \sigma_{\max} \frac{\omega_1 ( \overline{\fset \setminus t \mathcal{B}} )}{ \mu \sqrt{n}}, \notag
\end{equation}
where $\sigma_{\max}$ is defined as in Theorem \ref{thm_no_mismatch}.
\end{cor}

The proofs of Lemma \ref{lem_mismatched} and Corollary \ref{cor_mismatched} are similar to the proofs of Lemma \ref{lem_fundamental} and Theorem \ref{thm_no_mismatch}, respectively.

\section{Applications} \label{sec_appl}

Once the conditions (\ref{eq_RSC}) and (\ref{eq_RRSC}) are verified, our results Theorem \ref{thm_no_mismatch} and Corollary \ref{cor_mismatched} immediately follow. We explicitly verify the conditions for two applications and obtain the corresponding estimation error bounds.

The first application is regression by the constrained LS estimator in a Gaussian linear model. Let $\theta^\natural \in \mathbb{R}^p$ and $a_1, \ldots, a_n$ be vectors in $\mathbb{R}^p$. The sample is given by
\begin{equation}
y_i = \langle a_i, \theta^\natural \rangle + \sigma w_i, \quad i = 1, \ldots, n, \notag
\end{equation}
for some $\sigma > 0$, where $w_1, \ldots, w_n$ are i.i.d. standard Gaussian random variables. We consider the constrained LS estimator, for which $f_n$ is given by  (\ref{eq_LS}), and $\mathcal{G} := \set{ \theta: g ( \theta ) \leq c }$ for some $c \geq g ( \theta^\natural )$, where $g$ can be any convex continuous function.

\begin{cor} \label{cor_lasso_error}
Consider the Gaussian linear model and the constrained LS estimator described above. 
Assume that the entries of $a_1, \ldots, a_n$ are either all i.i.d. standard Gaussian or all i.i.d. Rademacher random variables.
Let $\epsilon \in ( 0, 1 )$. For any $t \geq 0$, there exist positive constants $c_1$ and $c_2$ such that if
\begin{equation}
\sqrt{n} \geq \frac{c_1 \alpha^2 \omega_1(\overline{\fset \setminus t \mathcal{B}})}{\epsilon}, \label{eq_lasso_sample_complexity}
\end{equation}
then we have
\begin{equation}
\mathbb{E}\, \norm{ \hat{\theta} - \theta^\natural }_2 \leq t + 2 \sqrt{ 2 \pi } \sigma\, \frac{\omega_1 ( \overline{\fset \setminus t \mathcal{B}} )}{ ( 1 - \epsilon ) \sqrt{n}} ,  \label{eq_lasso_error}
\end{equation}
with probability at least $1 - \exp(-c_2 \epsilon^2 n) > 1/2$ when $n$ is large enough.
\end{cor}

\begin{rem}
%When the constraint is matched, we can simply set $t = 0$ to minimize the error bound; one can check that $\ell_1 ( \overline{ \fset \setminus t \mathcal{B} } )$ does not change with $t$ if the constraint is matched. 
When the constraint is matched, we can simply set $t = 0$.
Recall that $t$ cannot be zero for the mismatched constraint case when $n < p$ (cf. Section \ref{sec_mismatch}). This remark also applies to Corollary \ref{cor_glm_error} below.
\end{rem}

\begin{rem}
For the mismatched constraint case, Corollary (\ref{cor_lasso_error}) is minimax optimal for the Lasso in the Gaussian linear model. We address this in Section \ref{sec_mismatch_further}.
\end{rem}

Corollary \ref{cor_lasso_error} is consistent with \cite{Oymak2013}. The result in \cite{Oymak2013} is sharper, while Corollary \ref{cor_lasso_error} is more general as it also covers the mismatched constraint case.

%When the constraint is matched, letting $t \rightarrow 0$ one may verify that the requirement on $n$ reduces to
%\begin{equation}
%\sqrt{n} \geq \frac{c \alpha^2 \ell_1 ( \fcone )}{\varepsilon} \notag 
%\end{equation}
%and the error bound (\ref{eq_lasso_error}) reduces to
%\begin{equation}
%\mathbb{E}\, \norm{ \hat{\theta} - \theta^\natural }_2 \leq \frac{1+\varepsilon}{(1-\varepsilon)^{3/2}}\sqrt{\frac{2\pi \omega_1(\fset)}{n}}
%\end{equation} which is consistent with existing results on constrained and regularized LS estimators (ignoring the constants). Note that our bound is more general as it also covers the mismatched constraint case.

The second application is $\ell_1$-constrained ML regression in a canonical GLM. %with an $\ell_1$-constraint. 

\begin{cor} \label{cor_glm_error}
Consider the canonical GLM and the constrained ML estimator described in Section \ref{sec_matched}, for $g ( \theta ) := \norm{ \theta }_1$ and $c \geq \norm{ \ttrue }$. Assume that $f_n$ in (\ref{eq_fn_glm}) is twice continuously differentiable, and the entries of $a_1, \ldots, a_n$ are i.i.d. Rademacher random variables. Let $\epsilon \in ( 0, 1 )$. For any $t \geq 0$, there exist positive constants $c_1$, and $c_2$ such that if (\ref{eq_lasso_sample_complexity}) is satisfied, then we have 
\begin{equation}
\mathbb{E}\, \norm{ \hat{\theta} - \theta^\natural }_2 \leq t +2 \sqrt{ 2 \pi } \, \sigma_{\max} \frac{\omega_1 ( \overline{\fset \setminus t \mathcal{B}} )}{ ( 1 - \epsilon ) \sqrt{n}}, \label{eq_glm_error}
\end{equation}
with probability at least $1 - \exp(c_2 \epsilon^2 n) > 1/2$ when $n$ is large enough, where $\sigma_{\max} := \max_i \sqrt{ \mathrm{var}\, y_i }$ is bounded above by a constant independent of $n$. 
\end{cor}

To the best of our knowledge, there are not existing results for $\ell_1$-constrained ML regression in GLMs. Here we compare Corollary \ref{cor_glm_error} with \cite{Negahban2010}, which provides an error bound for $\ell_1$-penalized ML estimators in GLMs
%\footnote{Recall that, however, there exist some fundamental differences between the constrained and penalized formulations (cf. Section \ref{sec_related_work}).}
. Recall that, however, the correspondence between the constrained and penalized estimators is currently unclear. When the constraint is matched and $\theta^\natural$ is $s$-sparse, Corollary \ref{cor_glm_error} states that when $n = \Omega ( s \log ( p / s ) )$,
\begin{equation}
\mathbb{E}\, \norm{ \hat{\theta} - \theta^\natural }_2 = O \left( \sqrt{\frac{s}{n} \log \left( \frac{p}{s} \right) } \right) \notag
\end{equation}
by Proposition 3.10 in \cite{Chandrasekaran2012}, which essentially coincides with Corollary 5 in \cite{Negahban2010}\footnote{We cite \cite{Negahban2010} instead of the published version \cite{Negahban2012}, because the estimation error bound only appears in \cite{Negahban2010}.}. We note that \cite{Negahban2010} only provides an error bound for the $\ell_1$-penalization case. 
%It is currently unclear to us whether extensions to other regularizers are possible or not. 

\section{Sharpness of Our Error Bound} \label{sec_mismatch_further}

It has been shown that in a Gaussian linear model with $\mathcal{G}$ being an $\ell_1$-ball, \emph{any} estimator $\hat{\theta}_{\text{arbitrary}}$ must satisfy, with probability larger than $1/2$, 
\begin{equation}
\max_{\theta^\natural \in \mathcal{G}} \norm{ \hat{\theta}_{\text{arbitrary}} - \theta^\natural }_2 = \Omega ( n^{-1/4} ), \notag
\end{equation}
under some technical conditions \cite{Raskutti2011}. Now we show our error bound for the Lasso in Corollary \ref{cor_lasso_error} actually achieves the error decaying rate $O ( n^{-1/4} )$ in the mismatched constraint case, and hence cannot be essentially improved.

By the definition of the Gaussian width, we have, for any $t > 0$,
\begin{equation}
\omega_1 \left( \overline{\fset \setminus t \mathcal{B}} \right) = \frac{ \omega_t \left( \overline{\fset \setminus t \mathcal{B}}  \right) }{t} = \frac{ \omega_t \left( \fset \right) }{t}, \notag
\end{equation} 
and hence the estimation error bound in Corollary \ref{eq_lasso_error} can be written as
\begin{equation}
\mathbb{E}\, \norm{ \that - \ttrue }_2 \leq t + \frac{C}{t} \frac{ \omega_t ( \fset ) }{\sqrt{n}}, \label{compare}
\end{equation}
for some $C > 0$, when $n$ is large enough such that (\ref{eq_lasso_sample_complexity}) is satisfied.

Define the \emph{global Gaussian width}: 
\begin{equation}
\omega ( \fset ) := \mathbb{E}\, \sup_{ v \in \fset } \set{ \left\langle h, x \right\rangle }, \notag
\end{equation}
where $h \in \mathbb{R}^p$ is a vector of i.i.d. standard Gaussian random variables. By definition, $\omega_t ( \fset )$ is bounded above by $\omega ( \fset )$, independent of $n$. Replacing $\omega_t ( \fset )$ by $\omega ( \fset )$ in (\ref{compare}), we have a looser error upper bound:
\begin{equation}
\mathbb{E}\, \norm{ \that - \ttrue }_2 \leq t + \frac{C}{t} \frac{\omega( \fset )}{\sqrt{n}}, \notag
\end{equation}
Minimizing this bound over all $t > 0$, we obtain the $O ( n^{-1/4} )$ error decaying rate. Similar discussion can be found in \cite{Plan2014a}.

%Of course, \cite{Negahban2010} does not need to consider the mismatched constraint case, due to the difference in problem formulations.

\section{Discussion}

Note that by the elementary argument in Section \ref{sec_framework}, we arrive at an estimation error bound (\ref{eq_concentration}) that holds \emph{surely}. It is possible to derive a concentration-type error guarantee based on this sure error bound, which we are working on.

Our framework is not restricted to constraint sets of the form (\ref{eq_that}); it applies to any non-empty closed convex set $\mathcal{G}$, as we only require $\iota_{\mathcal{G}} ( \cdot )$ to be proper closed convex in the proof. This observation is crucial to applying our framework to analyze constrained estimators for quantum tomography \cite{Flammia2012,Gross2010} and photon-limited imaging systems \cite{Raginsky2010}, which we are studying.

In this paper, we consider a random matrix $A$, and discuss the expected estimation error with respect to both $A$ and the sample $( y_1, \ldots, y_n )$. The extension to the the case where $A$ is deterministic is technically non-trivial, and we have not obtained a satisfactory result. We address this in the remark following the proof of Theorem \ref{thm_no_mismatch} in the appendix.

\section*{Acknowledgements}

This work was supported in part by the European Commission under Grant MIRG-268398, ERC Future Proof, SNF 200021-132548, SNF 200021-146750 and SNF CRSII2-147633.

\appendix

% Labels
\def \propRSCHessian {3.1}
\def \corMismatched {5.2}
\def \corLassoError {6.1}
\def \corGLMError {6.2}
\def \eqLassoError {16}
\def \eqRRSC {14}
\def \secMismatched {5}

\def \lemMismatched {5.1}
\def \thmNoMismatch {4.1}
\def \lemFundamental {3.2}

\section{Proof of Proposition {\propRSCHessian}}

We have
\begin{align}
\left\langle \nabla f_n ( \ttrue + e ) - \nabla f_n ( \ttrue ), e \right\rangle = \int_0^1 \left\langle e, \nabla^2 f_n ( \ttrue + \lambda e ) e \right\rangle \, d \lambda. \notag
\end{align}
The right-hand side is always larger than $\mu \norm{e}_2^2$ by assumption.

\section{Proof of Theorem \thmNoMismatch}
The main goal of the proof is to evaluate $\mathbb{E}\, \norm{ \Pi_{ \fcone } \left( - \nabla f_n ( \ttrue ) \right) }_2$. Here the expectation is with respect to both $A$ and the sample $( y_i )_{i = 1, \ldots, n}$.

We start with an equivalent formulation: 
\begin{equation}
\mathbb{E}\, \norm{ \Pi_{\fcone} \left( - \nabla f_n ( \ttrue ) \right) }_2 = \mathbb{E}\, \sup_{ v \in \fcone \cap \mathcal{S}^{p - 1} } \set{ \left\langle - \nabla f_n ( \ttrue ), v \right\rangle }, \label{eq_proj_sup}
\end{equation}
where $\mathcal{S}^{p - 1}$ denotes the unit $\ell_2$-sphere in $\mathbb{R}^p$. It is well known that in a canonical GLM, we have
\begin{equation}
\nabla f_n ( \ttrue ) = - \frac{1}{n} A^T \varepsilon, \label{eq_nablaF}
\end{equation}
where $\varepsilon := ( y_i - \mathbb{E}\, y_i )_{i = 1, \ldots, n}$, and hence
\begin{equation}
\mathbb{E}\, \norm{ \Pi_{\fcone} \left( - \nabla f_n ( \ttrue ) \right) }_2 = \frac{1}{n} \, \mathbb{E}\, \sup_{  v \in \fcone \cap \mathcal{S}^{p - 1}} \set{ \left\langle A^T \varepsilon,  v \right\rangle }. \notag
\end{equation}

To proceed, we need the following symmetrization inequality. The symmetrization inequality is different from the well-known symmetrization inequality by a Rademacher process, so we show it here for completeness.

\begin{lem}[\cite{Handel2014}] \label{lem_symmetrization}
Let $\xi_1, \ldots, \xi_n$ be independent real-valued random variables, and let $\mathcal{F}$ be a class of real functions. We have
\begin{equation}
\mathbb{E}\, \sup_{f \in \mathcal{F}} \set{ \sum_{i = 1}^n \left[ f ( \xi_i ) - \mathbb{E}\, f ( \xi_i ) \right] } \leq \sqrt{2 \pi}  \,\mathbb{E}\, \sup_{ f \in \mathcal{F} } \set{ \sum_{i = 1}^n h_i f ( \xi_i ) }, \notag
\end{equation}
where $h_1, \ldots, h_n$ are i.i.d. standard Gaussian random variables.
\end{lem}

\begin{rem}
In \cite{Handel2014}, the lemma is stated for the case when $\xi_1, \ldots, \xi_n$ are i.i.d. The case when $\xi_1, \ldots, \xi_n$ are not necessarily identical can be proved in a similar way, as noted in \cite{Pollard1984}.
\end{rem}

By Lemma \ref{lem_symmetrization}, we have
\begin{align}
\mathbb{E}\, \sup_{  v \in \fcone \cap \mathcal{S}^{p - 1}} \set{ \left\langle A^T \varepsilon, v \right\rangle } 
& \quad = \mathbb{E}\, \sup_{  v \in \fcone \cap \mathcal{S}^{p - 1}} \set{ \left\langle \varepsilon, A v \right\rangle } \notag \\
& \quad \leq \sqrt{2 \pi} \, \mathbb{E}\, \sup_{ v \in \fcone \cap \mathcal{S}^{p - 1} } \set{ \left\langle h \cdot \varepsilon, A v \right\rangle }, \notag
\end{align}
where $h \cdot \varepsilon := ( h_i \varepsilon_i )_{i = 1, \ldots, n}$, and $h_1, \ldots, h_n$ are i.i.d. standard Gaussian random variables. Note that $h \cdot \varepsilon$ is a random Gaussian vector with zero mean and covariance matrix $\Sigma \in \mathbb{R}^{n \times n}$ which is dependent on $A$ in general; moreover, since the entries in $\varepsilon$ are independent, $\Sigma$ is a diagonal matrix with diagonal entries given by $\Sigma_{i,i} := \mathrm{var}\, y_i$. Define $\tilde{h} := ( \tilde{h}_i )_{i = 1, \ldots, n}$, where $\tilde{h}_i := \Sigma_{i,i}^{-1/2} h_i \varepsilon_i$. Then $\tilde{h}$ is a vector of i.i.d. standard Gaussian random variables; furthermore, it is still a vector of i.i.d. standard Gaussian random variables condition on $A$, and hence it is statistically independent of $A$. 

Since $h \cdot \varepsilon$ and $\sqrt{\Sigma} \tilde{h}$ have the same probability distribution, we can write
\begin{equation}
\mathbb{E}\, \sup_{ v \in \fcone \cap \mathcal{S}^{p - 1} } \set{ \left\langle h \cdot \varepsilon, A v \right\rangle } = \mathbb{E}\, \sup_{v \in \fcone \cap \mathcal{S}^{p - 1}} \set{ \left\langle \sqrt{\Sigma} \tilde{h}, A v \right\rangle }. \notag
\end{equation}
Let $\mathcal{T} := \fcone \cap \mathcal{S}^{p - 1}$. Condition on any given $A$ (and hence $\Sigma$), we consider two mean-zero Gaussian processes $\set{ X_t }_{t \in \mathcal{T}}$ and $\set{ Y_t }_{t \in \mathcal{T}}$ defined as
\begin{equation}
X_t := \left\langle \sqrt{\Sigma} \tilde{h}, A t \right\rangle, \quad Y_t := \sigma_{\max} \left\langle \tilde{h}, A t \right\rangle, \notag
\end{equation}
where $\sigma_{max} := \max_i \Sigma_{i,i} = \max_i \sqrt{\mathrm{var}\, \varepsilon_i}$. We have, for any $t_1, t_2 \in \mathcal{T}$, 
\begin{equation}
\mathbb{E}\, \abs{ X_{t_1} - X_{t_2} }^2 = \norm{ \Sigma A ( t_1 - t_2 ) }_2^2 \leq \sigma_{\max}^2 \norm{ A ( t_1 - t_2 ) }_2^2 = \mathbb{E}\, \abs{ Y_{t_1} - Y_{t_2} }^2. \notag
\end{equation}
By Slepian's lemma, this implies 
\begin{equation}
\mathbb{E}\, \sup_{t \in \mathcal{T}} X_t \leq \mathbb{E}\, \sup_{t \in \mathcal{T}} Y_t. \notag
\end{equation}
Since the inequality holds given any realization of $A$, we have
\begin{align}
\mathbb{E}\, \sup_{  v \in \fcone \cap \mathcal{S}^{p - 1}} \set{ \left\langle A^T \varepsilon, v \right\rangle } 
&\leq \sqrt{2 \pi} \, \sigma_{\max} \, \mathbb{E}\, \sup_{v \in \fcone \cap \mathcal{S}^{p - 1}} \set{ \left\langle \tilde{h}, A v \right\rangle } \notag \\
&= \sqrt{2 \pi} \, \sigma_{\max} \, \mathbb{E}\, \sup_{v \in \fcone \cap \mathcal{S}^{p - 1}} \set{ \left\langle A^T \tilde{h}, v \right\rangle }. \notag
\end{align}

It remains to prove
\begin{equation}
\mathbb{E}\, \sup_{v \in \fcone \cap \mathcal{S}^{p - 1}} \set{ \left\langle A^T \tilde{h},  v \right\rangle } \leq \sqrt{n} \, \omega_1 ( \fcone ) := \sqrt{n} \, \mathbb{E}\, \sup_{ v \in \fcone \cap \mathcal{S}^{p - 1} }\set{ \left\langle \tilde{h}, v \right\rangle }. \label{eq_key}
\end{equation}
We consider two cases: 

\paragraph{Case 1: } If $A$ has i.i.d. standard Gaussian entries, then condition on $\tilde{h}$, $A^T \tilde{h}$ is a vector of mean-zero Gaussian random variables with covariance matrix $\norm{ \tilde{h} }_2 I$, and hence has the same probablity distribution as $\norm{ \tilde{h} } \bar{h}$, where $\bar{h}$ is a vector of i.i.d. standard Gaussian random variables independent of $\tilde{h}$. Therefore,
\begin{align}
\mathbb{E}\, \sup_{v \in \fcone \cap \mathcal{S}^{p - 1}} \set{ \left\langle A^T \tilde{h},  v \right\rangle } &= \mathbb{E} \, \sup_{v \in \fcone \cap \mathcal{S}^{p - 1}} \set{ \left\langle \norm{ \tilde{h} } \bar{h} , v \right\rangle } \notag \\
&= \left( \mathbb{E}_{\tilde{h}}\, \norm{ \tilde{h} }_2 \right) \, \mathbb{E}_{\bar{h}} \sup_{v \in \fcone \cap \mathcal{S}^{p - 1}} \set{ \left\langle \bar{h}, v \right\rangle } \notag \\
&\leq \sqrt{n}\, \omega_1 ( \fcone ). \notag
\end{align}

\paragraph{Case 2: } If $A$ has i.i.d. Rademacher entries, then condition on $A$, $A^T \tilde{h}$ is a vector of mean-zero Gaussian random variables with covariance matrix $n I$, and hence has the same probability distribution as $\sqrt{n} \bar{h}$, where $\bar{h}$ is a vector of i.i.d. standard Gaussian random variables. Therefore, 
\begin{align}
\mathbb{E}\, \sup_{v \in \fcone \cap \mathcal{S}^{p - 1}} \set{ \left\langle A^T \tilde{h},  v \right\rangle } &= \mathbb{E} \, \sup_{v \in \fcone \cap \mathcal{S}^{p - 1}} \set{ \left\langle \sqrt{n} \bar{h} , v \right\rangle } \notag \\
&= \sqrt{n}\, \omega_1 ( \fcone ). \notag
\end{align}

In summary, we obtain
\begin{equation}
\mathbb{E}\, \norm{ \Pi_{\fcone} ( - \nabla f_n ( \ttrue ) ) }_2 \leq \sqrt{ 2 \pi }\, \sigma_{\max}\, \frac{ \omega_1 ( \fcone ) }{ \sqrt{n} }, \notag
\end{equation}
if the entries of $A$ are i.i.d. standard Gaussian or Rademacher random variables, for a canonical GLM, where the expectation is with respect to both $A$ and the sample $( y_i )_{i = 1, \ldots, n}$.

Let $\mathcal{E}$ denote that event that the RSC condition holds. Then we have
\begin{align}
\mathbb{E}\, \norm{ \Pi_{\fcone} ( - \nabla f_n ( \ttrue ) ) }_2 = & \, \mathbb{P} ( \mathcal{E} ) \, \mathbb{E}_{A, ( y_i ) \vert \mathcal{E}}\, \norm{ \Pi_{\fcone} ( - \nabla f_n ( \ttrue ) ) }_2 \notag \\
& \quad + \mathbb{P} ( \mathcal{E}^C ) \, \mathbb{E}_{A, ( y_i ) \vert \mathcal{E}^C}\, \norm{ \Pi_{\fcone} ( - \nabla f_n ( \ttrue ) ) }_2, \notag
\end{align}
and hence
\begin{align}
\mathbb{E}_{A, ( y_i ) \vert \mathcal{E}}\, \norm{ \Pi_{\fcone} ( - \nabla f_n ( \ttrue ) ) }_2 &\leq \frac{ \mathbb{E}\, \norm{ \Pi_{\fcone} ( - \nabla f_n ( \ttrue ) ) }_2 }{ \mathbb{P} ( \mathcal{E} ) } \notag \\
&\leq 2 \mathbb{E}\, \norm{ \Pi_{\fcone} ( - \nabla f_n ( \ttrue ) ) }_2, \notag
\end{align}
where we applied the assumption that $\mathbb{P} ( \mathcal{E} ) \geq 1/2$. By Lemma \lemFundamental, this implies 
\begin{align}
\mathbb{E}_{A, \varepsilon \vert \mathcal{E}}\, \norm{ \that - \ttrue }_2 
&\leq \frac{1}{\mu} \mathbb{E}_{A, ( y_i ) \vert \mathcal{E}}\, \norm{ \Pi_{\fcone} ( - \nabla f_n ( \ttrue ) ) }_2 \notag \\
&\leq 2 \sqrt{2 \pi} \, \sigma_{\max} \, \frac{ \omega_1 ( \fcone ) }{\mu \sqrt{n}}. \notag
\end{align}
This completes the proof.

\begin{rem}
If we want to adapt this proof to the deterministic $A$ case, a technical issue arises when bounding the right-hand side of (\ref{eq_key}). As the random process $\set{ \tilde{X}_v := \left\langle \tilde{h}, v \right\rangle }_{v \in \mathcal{V}}$, where $\mathcal{V} := \fcone \cap \mathcal{S}^{p - 1}$, is a mean-zero Gaussian process, a standard approach is to bound $\sup_{v \in \mathcal{V}} \tilde{X}_v$ by Slepian's lemma. Note that, for any $v_1, v_2 \in \mathcal{V}$,
\begin{equation}
\mathbb{E}\, \abs{ \tilde{X}_{v_1} - \tilde{X}_{v_2} }^2 = \norm{ A ( v_1 - v_2 ) }_2^2, \notag
\end{equation}
and hence an upper-bound on $\mathbb{E}\, \abs{ \tilde{X}_{v_1} - \tilde{X}_{v_2} }^2$ would depend on the largest eigenvalue of $A$. The largest eigenvalue of $A$, however, cannot be bounded above by a constant independent of $n$ under the high-dimensional setting. Although we can weaken the requirement on $A$ to a restricted smoothness condition as
\begin{equation}
\norm{ A v }_2 \leq \sqrt{ 1 + \epsilon } \norm{v}_2, \quad \text{for all } v \in \fcone \cap \mathcal{S}^{p - 1}, \notag
\end{equation}
which, by Theorem \ref{thm_mendelson}, holds with high probability. This condition does not imply 
\begin{equation}
\norm{ A ( v_1 - v_2 ) }_2^2 \leq C \norm{ v_1 - v_2 }_2^2, \notag
\end{equation}
for some dimension-independent constant $C > 0$, for all $v_1, v_2 \in \mathcal{V}$.
\end{rem}

\section{Proof of Lemma \lemMismatched} \label{sec_proof_lem_mismatch}

Let $e := \that - \ttrue$. If $e \in \fset \setminus t \mathcal{B}$, following the proof of Theorem \thmNoMismatch, we obtain
\begin{equation}
\norm{e}_2 \leq \frac{1}{\mu} \norm{ \Pi_{\overline{\fset \setminus t \mathcal{B}}} \left( - \nabla f_n ( \ttrue ) \right) }_2, \notag
\end{equation}
where $\overline{\fset \setminus t \mathcal{B}}$ denotes the conic hull of $\fset \setminus t \mathcal{B}$. If $e \in t \mathcal{B}$, we have the na\"{\i}ve bound: $\norm{e}_2 \leq t$. Therefore, 
\begin{align}
\norm{e}_2 &\leq \max \set{ t, \frac{1}{\mu} \norm{ \Pi_{\overline{\fset \setminus t \mathcal{B}}} \left( - \nabla f_n ( \ttrue ) \right) }_2 } \notag \\
&\leq t + \frac{1}{\mu} \norm{ \Pi_{\overline{\fset \setminus t \mathcal{B}}} \left( - \nabla f_n ( \ttrue ) \right) }_2. \notag
\end{align}
The lemma follows by taking expectations on both sides.

\section{Proof of Corollary \corMismatched} \label{sec_proof_cor_mismatch}

Let $e := \that - \ttrue$. If $e \in \fset \setminus t \mathcal{B}$, following the proof of Theorem \thmNoMismatch, we can obtain
\begin{equation}
\mathbb{E}\, \norm{ e }_2 \leq 2 \sqrt{ 2 \pi } \, \sigma_{\max} \frac{\omega_1 ( \overline{ \fset \setminus t \mathcal{B} } )}{ \mu \sqrt{n}}; \notag
\end{equation}
otherwise, we can bound the expected estimation error from above by $t$. Therefore, 
\begin{align}
\mathbb{E}\, \norm{e}_2 &\leq \max \set{ t, 2 \sqrt{ 2 \pi } \, \sigma_{\max} \frac{\omega_1 ( \overline{ \fset \setminus t \mathcal{B} } )}{ \mu \sqrt{n}} } \notag \\
& \leq t + 2 \sqrt{ 2 \pi } \, \sigma_{\max} \frac{\omega_1 ( \overline{ \fset \setminus t \mathcal{B} } )}{ \mu \sqrt{n}}. \notag
\end{align}

\section{Proof of Corollary \corLassoError \ and Corollary \corGLMError}

%%%%%%%%%
%%%%%%%%%

The proofs in this section rely on the following theorem \cite{Mendelson2007}.

\begin{thm}[\cite{Mendelson2007}] \label{thm_mendelson}
Let $\mathcal{T} \subseteq \mathbb{R}^p$ be star-shaped. Let $A \in \mathbb{R}^{n \times p}$, $n < p$, whose rows are i.i.d. isotropic subgaussian random vectors with subgaussian norm $\alpha \geq 1$, and let $\epsilon \in ( 0, 1 )$. Then there exist constants $c_1$ and $c_2$ such that for all $x \in \mathcal{T}$ satisfying
\begin{equation}
\norm{ x }_2 \geq \gamma_n^* \left( \frac{\epsilon}{c_1 \alpha^2}, \mathcal{T} \right) := \inf \set{ t > 0: t \geq \frac{c_1 \alpha^2 \omega_t ( \mathcal{T} )}{\epsilon \sqrt{n}} }, \label{eq_rkstar}
\end{equation}
we have
\begin{equation}
( 1 - \epsilon ) \norm{ x }_2^2 \leq \frac{\norm{ A x }_2^2}{n} \leq ( 1 + \epsilon ) \norm{ x }_2^2 \notag
\end{equation}
with probability at least $1 - \exp \left( - c_2 \epsilon^2 n / \alpha^4 \right)$.
\end{thm}

We note that the sub-Gaussian norm of a vector of i.i.d. standard Gaussian entries or i.i.d. Rademacher entries is bounded above by a constant \cite{Vershynin2012}.

%Recall that $\ell_t ( \cdot )$ is the Gaussian width defined in Definition \ref{def_gwidth}. Examples of an isotropic subgaussian random vector include a vector of i.i.d. standard Gaussian random variables, and a vector of i.i.d. Rademacher random variables \cite{Vershynin2012}. 

\subsection{Proof of Corollary \corLassoError} \label{sec_Lasso}

We prove by Corollary \corMismatched.

Let $A$ be defined as in Theorem \thmNoMismatch. We verify the condition (\eqRRSC) by Theorem \ref{thm_mendelson}. Since $\omega_t ( \overline{\fset \setminus t \mathcal{B}} ) = t \omega_1 ( \overline{\fset \setminus t \mathcal{B}} )$, the condition (\ref{eq_rkstar}) is equivalent to requiring 
\begin{equation}
\sqrt{n} \geq \frac{c_1 \alpha^2 \omega_1 ( \overline{\fset \setminus t \mathcal{B}} )}{\epsilon}. \notag
\end{equation}
Once this inequality is satisfied, we can set $\mu = 1 - \epsilon$, and the condition (\eqRRSC) hold with probability at least $1 - \exp \left( - c_2 \epsilon^2 n / \alpha^4 \right)$. Note that $\sigma_{\max} = \sqrt{ \mathbb{E}\, w_i^2 } = \sigma$. This completes the proof.

\subsection{Proof of Corollary \corGLMError}

We prove the corollary by Corollary \corMismatched.

It is known that
\begin{equation}
\nabla^2 f_n ( \theta ) = \frac{1}{n} A^T D(\theta) A \notag
\end{equation}
for the ML estimator in a canonical GLM, where $A$ is defined as in Theorem \thmNoMismatch, and $D( \theta )$ is a diagonal matrix; furthermore, there exists a continuous strictly positive function $\phi$ such that the $( i, i )$-th entry  of $D ( \theta )$ is given by $\phi ( \langle a_i, \theta \rangle )$. Since the entries of $A$ are i.i.d. Rademacher random variables, for any $\theta \in \mathcal{G}$, 
\begin{equation}
\abs{ \left\langle a_i, \theta \right\rangle } \leq \norm{ a_i }_{\infty} \norm{ \theta }_1 \leq c. \notag
\end{equation}
By the extreme value theorem, the diagonal entries of $D(\theta)$ are bounded below by a constant $\nu > 0$ for all $\theta \in \mathcal{G}$, which is independent of $n$. Similarly, $\sigma_{\max}$ is bounded above by a constant independent of $n$. 

The rest of the proof is similar to the last paragraph in the previous sub-section. 
By Theorem \ref{thm_mendelson}, if we choose $n$ such that 
\begin{equation}
\sqrt{n} \geq \frac{c \alpha^2 \omega_1 ( \overline{\fset \setminus t \mathcal{B}} )}{\epsilon}, \notag
\end{equation}
then the condition (\eqRRSC) holds with probability at least $1 - \exp \left( - c_2 \epsilon^2 n / \alpha^4 \right)$ with $\mu = \nu ( 1 - \epsilon )$.

{\small
\bibliography{list}
\bibliographystyle{IEEEtranS}
}

\end{document}